\documentclass[aos,MSNbibl,nameyear,dvips]{arximspdf}

\doi{10.1214/13-AOS1148} 
\volume{41}
\issue{4}
\pubyear{2013}
\firstpage{2218}
\lastpage{2235}

\makeatletter
\def\cal{\mathcal}

\def\open{}
\def\E{{\open\mathbb{E}}}
\def\I{{\open\mathbb{I}}}
\def\proj{\operatorname{pr}^{\bot}}

\newtheorem{theorem}{Theorem}
\newproclaim{remark}{Remark}
\newtheorem{corol}{Corollary}
\newtheorem{proposition}{Proposition}
\makeatother

\begin{document}
\begin{frontmatter}

\title{Optimal crossover designs for the proportional model}
\runtitle{Design for proportional model}

\begin{aug}
\author{\fnms{Wei} \snm{Zheng}\corref{}\ead[label=e1]{weizheng@iupui.edu}}
\runauthor{W. Zheng}
\affiliation{Indiana University-Purdue University Indianapolis}
\address{Department of Mathematical Sciences\\
Indiana University-Purdue University Indianapolis\\
Indianapolis, Indiana 46202-3216\\
USA\\
\printead{e1}}
\end{aug}

\received{\smonth{4} \syear{2013}}
\revised{\smonth{7} \syear{2013}}


\begin{abstract}
In crossover design experiments, the proportional model, where the
carryover effects are proportional to their direct treatment effects,
has draw attentions in recent years. We discover that the universally
optimal design under the traditional model is E-optimal design under
the proportional model. Moreover, we establish equivalence theorems of
Kiefer--Wolfowitz's type for four popular optimality criteria, namely A,
D, E and T (trace).
\end{abstract}

%

\end{frontmatter}

\section{Introduction}

Let $\Omega_{p,t,n}$ be the collection of all crossover designs with
$p$ periods, $t$ treatments, and $n$ subjects. In an experiment based
on design $d\in\Omega_{p,t,n}$, the response from subject $u\in\{
1,2,\ldots,n\}$ in period $k\in\{ 1,2,\ldots,p\}$, to which treatment
$d(k,u)\in\{ 1,2,\ldots,t\}$ was assigned by design $d$, is traditionally
modeled as
%
%
\begin{equation}
\label{model1} Y_{dku}=\mu+\alpha_k+\beta_u+
\tau_{d(k,u)}+\gamma _{d(k-1,u)}+\varepsilon_{ku}.
\end{equation}
Here, $\mu$ is the general mean,
$\alpha_k$ is the $k$th period effect, $\beta_u$ is the $u$th subject
effect, $\tau_{d(k,u)}$ is the (direct) effect of treatment $d(k,u)$,
and $\gamma_{d(k-1,u)}$ is the
carryover effect of treatment $d(k-1,u)$ that subject $u$ received in
the previous period (by convention $\gamma_{d(0,u)}=0$). A central
problem in the area of crossover design is to find the best design
among $\Omega_{p,t,n}$ for estimating the direct, and sometimes also
carryover, treatment effects. Since Hedayat and Afsarinejad (\citeyear{HA75}, \citeyear{HedAfs78}) the optimal design problems have been mainly studied
under model
(\ref{model1}). Examples include Cheng and Wu (\citeyear{CheWu80}), \citet{Kun84},
\citet{Stu91}, Hedayat and Yang (\citeyear{HedYan03}, \citeyear{HedYan04}) and \citet{HedZhe10} among others. For approximate design solutions,
see Kushner
(\citeyear{Kus97,Kus98}), \citet{KunMar00}, \citet{KunStu02}, and
\citet{Zh2013} among others.

Many variants of model (\ref{model1}) have been proposed in literature.
The main focus is on different modelings of carryover effects, such as
no carryover effects model, mixed carryover effects model [\citet{KunStu02}] and the full interaction model [\citet{Paretal11}]. The
choice of model should be based on practical background and it is the
responsibility of design theorists to provide recipes of optimal or
efficient designs for each of these models. Here we consider model
(\ref
{model2}) below because usually (i) it is essential to choose a
parsimonious but reasonable model; (ii) The treatment having the
larger direct effect in magnitude usually yields the larger carryover
effect:
%
%
\begin{equation}
\label{model2} Y_{dku}=\mu+\alpha_k+\beta_u+
\tau_{d(k,u)}+\lambda\tau _{d(k-1,u)}+\varepsilon_{ku}.
\end{equation}

Throughout the paper, we call this model as the \textit{proportional}
model. \citet{KemFerDav01} proposed this model and some
theoretical results are later derived by \citet{BaiKun06} and
\citet{BosStu07}. The main difficulty is due to the nonlinear
term $\lambda\tau_{d(k-1,u)}$ in the model. In this paper, we show that
universally optimal designs for estimating treatment effects under the
traditional model is E-optimal under the proportional model regardless
the value of $\lambda$. Unlike the traditional model, the proportional
model do not yield universally optimal designs in general. Instead, we
derive equivalence theorems for four popular optimality criteria,
namely A, D, E and T. Besides, we derive optimal designs for estimating
$\lambda$.

The rest of this paper is organized as follows. Section~\ref{sec:2}
briefly introduces the universal optimality of Kushner's design under
the traditional linear model as well as some necessary notation to be
used for the rest of this paper. Section~\ref{sec:3} studies the
optimal design problem for the proportional model. Finally, Section~\ref{sec:exm}
gives some examples of optimal designs under different values
of $p$ and $t$.

\section{Some notation and Kushner's design}\label{sec:2}
Let $G$ be a temporary object whose meaning differs from context to
context. For a square matrix $G$, we define $G'$, $G^-$ and $\operatorname{tr}(G)$ to
represent the transpose, $g$-inverse and trace of $G$, respectively. The
projection operator $\proj$ is defined as $\proj G=I-G(G'G)^-G'$. For
two square matrices of equal size, $G_1$ and $G_2$, $G_1\leq G_2$ means
that $G_2-G_1$ is nonnegative definite. For a set $G$, the number of
elements in the set is represented by $|G|$. Besides, $I_k$ is the
$k\times k$ identity matrix and $1_k$ is the vector of length $k$ with
all its entries as $1$. We further define $J_k=1_k1_k'$ and
$B_k=I_k-J_k/k$. Finally, $\otimes$ represents the Kronecker product of
two matrices.

Let $Y_d=(Y_{d11},Y_{d21},\ldots,Y_{dp1},Y_{d12},\ldots,Y_{dpn})'$ be the
$np\times1$ response vector, then model (\ref{model1}) has the matrix
form
%
%
\begin{equation}
\label{model3} Y_d={ 1}_{np}\mu+Z{\alpha}+U{\beta}+{
T}_d{\tau}+{ F}_d{\gamma }+\varepsilon,
\end{equation}
where ${\alpha}=(\alpha_1,\ldots,\alpha_p)'$, ${\beta}=(\beta
_1,\ldots,\beta
_n)'$, ${\tau}=(\tau_1,\ldots,\tau_t)'$, ${\gamma}=(\rho_1,\ldots,\rho_t)'$,
$Z={ 1}_n\otimes{ I}_p$, $U=I_n\otimes1_p$, and ${ T}_d$ and ${ F}_d$
denote the treatment/subject and carryover/subject incidence matrices.
Here we assume $\E(\varepsilon)=0$ and $ \operatorname{Var}(\varepsilon
)=I_{n}\otimes
\Sigma$, where $\Sigma$ is a nonsingular within subject covariance
matrix. Define $\Sigma^{-1/2}$ to be the matrix such that $\Sigma
^{-1}=\Sigma^{-1/2}\Sigma^{-1/2}$. Let $\tilde{T}_d=I_n\otimes
\Sigma
^{-1/2}T_d$, $\tilde{F}_d=I_n\otimes\Sigma^{-1/2}F_d$, $\tilde
{Z}=I_n\otimes\Sigma^{-1/2}Z$ and $\tilde{U}=I_n\otimes\Sigma
^{-1/2}U$. The information matrix for the direct treatment effect $\tau
$ under model (\ref{model3}) is
\begin{eqnarray*}
C_d&=&\tilde{T}_d'\proj\bigl(\tilde{Z}|
\tilde{U}|\tilde{F}_d\bigr)\tilde {T}_d
\\
&=&C_{d11}-C_{d12}C_{d22}^-C_{d21},
\end{eqnarray*}
where $C_{dij}=G_i'\proj(\tilde{Z}|\tilde{U})G_j, 1\leq i,j\leq2$
with $G_1=\tilde{T}_d$ and $G_2=\tilde{F}_d$. Define $\tilde
{B}=\Sigma
^{-1}-\Sigma^{-1}J_p\Sigma^{-1}/1_p'\Sigma^{-1}1_p$, and note that
$\Sigma=I_p$ implies $\tilde{B}=B_p$. Straightforward calculations show
that $C_{dij}=G_i'(B_n\otimes\tilde{B})G_j, 1\leq i,j\leq2$ with
$G_1=T_d$ and $G_2=F_d$. A design is said to be universally optimal
[Kiefer (\citeyear{KI75})] if it maximizes $\Phi(C_d)$ for any $\Phi$ satisfying:
\begin{longlist}[(C.1)]
\item[(C.1)] $\Phi$ is concave.
\item[(C.2)] $\Phi(S'CS)=\Phi(C)$ for any permutation matrix $S$.
\item[(C.3)] $\Phi(bC)$ is nondecreasing in the scalar $b>0$.
\end{longlist}

In approximate design theory, a design $d\in\Omega_{p,t,n}$ is
considered as the result of selecting $n$ elements with replacement
from ${\cal S}$, the collection of all possible $t^p$ treatment
sequences. Now define the treatment sequence proportion $p_s=n_s/n$,
where $n_s$ is the number of replications of sequence $s$ in the
design. A~design in approximate design theory is then identified by the
vector $P_d=(p_s,s\in{\cal S})$ with the restrictions of $\sum_{s\in
{\cal S}}p_s=1$ and $p_s\geq0$.

Let $T_s$ (resp., $F_s$) be the $p\times t$ matrix $T_d$ (resp., $F_d$)
when $d$ consists of a single sequence $s$. For sequence $s\in{\cal
S}$ define $\hat{C}_{sij}=B_tG_i'\tilde{B}G_jB_t, 1\leq i,j\leq2$ with
$G_1=T_s$ and $G_2=F_s$. By direct calculations, we have
%
%
\begin{equation}
C_{dij}=\hat{C}_{dij}-nG_i'
\tilde{B}G_j,\qquad 1\leq i,j\leq2,
\end{equation}
where $\hat{C}_{dij}=n\sum_{s\in{\cal S}}p_s\hat{C}_{sij},1\leq
i,j\leq2$ with $G_1=\sum_{s\in{\cal S}}p_sT_sB_t$ and $G_2=\sum_{s\in
{\cal S}}p_sF_sB_t$. Further, we define $c_{sij}=\operatorname{tr}(\hat{C}_{sij})$,
$c_{dij}=\operatorname{tr}(\hat{C}_{dij})=n\sum_{s\in{\cal S}}p_sc_{sij}$, the
quadratic function $q_s(x)=c_{s11}+2c_{s12}x+c_{s22}x^2$, $q(x)=\max_{s}q_s(x)$, $y^*=\min_{-\infty<x<\infty}q(x)$, $x^*$ to be the unique
solution of $q(x)=y^*$ and ${\cal Q}=\{s\in{\cal S}|q_s(x^*)=y^*\}$.
\citet{Kus97} derived the following theorem.
%
\begin{theorem}[{[\citet{Kus97}]}]\label{thm:1222}
A design $d$ is universally optimal under model~(\ref{model3}) if and
only if
%
%
\begin{eqnarray}
\sum_{s\in{\cal Q}}p_s\bigl[
\hat{C}_{s11}+x^*\hat{C}_{s12}\bigr]&=&\frac
{y^*}{t-1}B_t,
\label{eqn:1231}
\\
\sum_{s\in{\cal Q}}p_s\bigl[
\hat{C}_{s21}+x^*\hat{C}_{s22}\bigr]&=&0,
\\
\sum_{s\in{\cal Q}}p_s\tilde{B}
\bigl(T_s+x^*F_s\bigr)B_t&=&0,
\\
\sum_{s\in{\cal Q}}p_s&=&1,
\\
p_s&=&0\qquad\mbox{if }s\notin{\cal Q}.\label{eqn:12312}
\end{eqnarray}
\end{theorem}

Let $\sigma$ be a permutation of the symbols $\{1,2,\ldots,t\}$. For a
sequence $s=(t_1,\ldots,t_p)$, we define $\sigma s=(\sigma
(t_1),\ldots,\sigma
(t_p))$. Note that $q_s(x)$ is invariant to treatment permutations,
that is,
%
%
\begin{equation}
\label{eqn:2042} q_s(x)=q_{\sigma s}(x),\qquad \sigma\in{\cal P}.
\end{equation}
Define the design $\sigma d$ by $P_{\sigma d}=(p_{\sigma^{-1}s},s\in
{\cal S})$. A design $d$ is said to be \textit{symmetric} if
$P_d=P_{\sigma d}$. Also we define symmetric blocks as $\langle
s\rangle
=\{\sigma s, \sigma\in{\cal P}\}$ where ${\cal P}$ is the collection
of all possible $t!$ permutations. For a symmetric design, we have
$p_{\tilde{s}}=p_{\langle s\rangle}/|\langle s\rangle|$ for any
$\tilde
{s}\in\langle s\rangle$, where $p_{\langle s\rangle}=\sum_{\tilde
{s}\in\langle s\rangle}p_{\tilde{s}}$. Given $p,t,n$, a symmetric
design $d$ is uniquely determined by $(p_{\langle s\rangle},s\langle
\in
\rangle{\cal S})$, where $s\langle\in\rangle{\cal S}$ means that $s$
runs through all distinct symmetric blocks contained in ${\cal S}$.
Equation (\ref{eqn:2042}) is essential for the following theorem.

\begin{theorem}[{[\citet{Kus97}]}]\label{thm:12222}
A symmetric design is universally optimal under model (\ref{model3}) if
%
%
\begin{eqnarray}
\sum_{s\in{\cal Q}}p_sq_s'
\bigl(x^*\bigr)&=&0,\label{eqn:12222}
\\
\sum_{s\in{\cal Q}}p_s&=&1,
\\
p_s&=&0,\qquad\mbox{if }s\notin{\cal Q},\label{eqn:12223}
\end{eqnarray}
where $q_s'(x)$ is the derivative of $q_s(x)$ with respective to $x$.
\end{theorem}

\section{Proportional model}\label{sec:3}

\subsection{Problem formulation and literature review}\label{sec:3.1}
We are interested in model~(\ref{model2}), which could be rewritten in
the matrix form
%
%
\begin{equation}
\label{model4} Y_d={ 1}_{np}\mu+T_d\tau+
\lambda F_d\tau+Z\alpha+U\beta +\varepsilon.
\end{equation}
Here we assume $\varepsilon\sim N(0,I_{n}\otimes\Sigma)$.
Fisher's information matrix for $\tau$ is
%
%
\begin{eqnarray}
\label{eqn:1204} C_{d,\tau_0,\lambda_0}(\tau)&=&(\tilde{T}_d+
\lambda_0 \tilde {F}_d)'\proj\bigl(\tilde{Z}|
\tilde{U}|\tilde{F}_d\tau_0\bigr) (\tilde {T}_d+
\lambda _0 \tilde{F}_d)
\nonumber
\\
&=&C_{d11}+\lambda_0(C_{d12}+C_{d21})+
\lambda_0^2C_{d22}
\\
&&-(C_{d12}+\lambda_0C_{d22})\tau_0
\bigl(\tau_0'C_{d22}\tau_0
\bigr)^{-1}\tau _0'(C_{d21}+
\lambda_0C_{d22}).\nonumber
\end{eqnarray}

Unlike model (\ref{model3}), model (\ref{model4}) is nonlinear, and
therefore the choice of optimal designs depends on the unknown
parameters $\lambda_0$ and $\tau_0$; see (\ref{eqn:1204}). The
nonlinearity of the model imposes the major difficulty on the problem.
For this, \citet{BosStu07} assumes that $\lambda_0$ is a known
parameter at the stage of data analysis, in which case the (Fisher's)
information matrix does not depend on $\tau_0$ and hence the same for
the choice of optimal designs. But such strategy inevitably yields
significant bias in the analysis stage when one do not have sufficient
knowledge about~$\lambda_0$.

Note that \citet{KemFerDav01} and \citet{BaiKun06}
also worked on $C_{d,\tau_0,\lambda_0}(\tau)$ even though they derived
it from the aspect of model approximation [\citet{FedHac97}, page~18] without normality assumption. For unknown $\tau_0$ and $\lambda_0$,
they adopted the following Bayesian type of criteria:
%
%
\begin{eqnarray}
\label{eqn:0103} \phi_{g,\lambda_0}(d)&=&\int\Phi\bigl(C_{d,\tau_0,\lambda_0}(\tau )
\bigr)g(\tau _0)\,d(\tau_0)
\nonumber
\\[-8pt]
\\[-8pt]
\nonumber
&=&\E_g\bigl(\Phi\bigl(C_{d,\tau_0,\lambda_0}(\tau)\bigr)\bigr),
\end{eqnarray}
where $g$ is the prior distribution of $\tau_0$. Note that they only
considered the special case of $\Sigma=I_p$ and $\Phi$ being the
A-criterion function. Particularly, \citet{KemFerDav01}
gave a search algorithm for A-efficient designs when $g$ is the density
function of a special multivariate normal distribution. \citet{BaiKun06}
proved the optimality of totally balanced design [\citet{KunStu02}] when $\Sigma=I_p$, $3\leq p\leq t$, the distribution
$g$ is exchangeable, and $-1\leq\lambda_0<\lambda^*$ with
%
%
\begin{equation}
\label{eqn:1223} \lambda^*=\frac{1}{p-1}-\frac{pt-t-1}{(p-1)(t-2)(pt-t-1-t/p)^2}.
\end{equation}

Note that $0<\lambda^*<1/(p-1)$. Hence the results of \citet{BaiKun06} will not be applicable when $p\geq t$ or the carryover effects
is positively proportional to the direct treatment effects with a
moderate or even larger magnitude. Here, we develop tools for finding
optimal designs for any value of $\lambda_0$ and $\Sigma$ and for four
popular criteria, namely A, D, E and T. For E-criterion, the optimal
design does not depend on the value of $\lambda_0$.

\subsection{Preliminary results}\label{sec:3.2}

Recall that the design $\sigma d$ is defined by $P_{\sigma
d}=(p_{\sigma
^{-1}s},s\in{\cal S})$. Let $S_{\sigma}$ be the unique permutation
matrix such that $T_{\sigma d}=T_dS_{\sigma}$ and $F_{\sigma
d}=F_dS_{\sigma}$ for any design $d$. Also define
\[
\sigma\tau_0=S_{\sigma}\tau_0.
\]
Let $\delta_{\tau_0}$ be the probability measure which puts equal mass
to each element in $\{\sigma\tau_0|\sigma\in{\cal P}\}$. We shall
focus on the special case of $g=\delta_{\tau_0}$ and then extend the
results to any arbitrary exchangeable distribution $g$. By definition
we have
\[
\phi_{\delta_{\tau_0},\lambda_0}(d)=\frac{1}{t!}\sum_{\sigma
}
\Phi \bigl(C_{d,\sigma\tau_0,\lambda_0}(\tau)\bigr),
\]
where the summation runs through all $t!$ permutations. Now we have:

\begin{theorem}\label{thm:1121}
In approximate design theory, given any values of the real number
$\lambda_0$ and the vector $\tau_0$, for any design $d$ there exists a
symmetric design, say $d^*$, such that
%
%
\begin{equation}
\label{eqn:12313} \phi_{\delta_{\tau_0},\lambda_0}(d)\leq\phi_{\delta_{\tau
_0},\lambda_0}\bigl(d^*\bigr).
\end{equation}
\end{theorem}
\begin{pf}
First, we observe that
%
%
\begin{equation}
\label{eqn:1121} C_{\sigma d,\tau_0,\lambda_0}(\tau)=S_{\sigma}'C_{ d,\sigma\tau
_0,\lambda_0}(
\tau)S_{\sigma}.
\end{equation}
For any given permutation $\sigma_0$, by (\ref{eqn:1121}) we have
%
%
\begin{eqnarray}
\label{eqn:11214} \phi_{\delta_{\tau_0},\lambda_0}(\sigma_0 d)&=&
\frac{1}{t!}\sum_{\sigma
}\Phi\bigl(C_{\sigma_0d, \sigma\tau_0,\lambda_0}(
\tau)\bigr)
\nonumber
\\
&=&\frac{1}{t!}\sum_{\sigma}\Phi
\bigl(S_{\sigma_0}'C_{d,\sigma_0
\sigma\tau
_0,\lambda_0}(\tau)S_{\sigma_0}\bigr)
\nonumber
\\[-8pt]
\\[-8pt]
\nonumber
&=&\frac{1}{t!}\sum_{\sigma}\Phi
\bigl(C_{d,\sigma_0 \sigma\tau
_0,\lambda
_0}(\tau)\bigr)
\\
&=&\phi_{\delta_{\tau_0},\lambda_0}(d).\nonumber
\end{eqnarray}
By direct calculations, we have
%
%
\begin{eqnarray}
\label{eqn:1119} C_{d,\tau_0,\lambda_0}(\tau,\lambda,\alpha)&=&\bigl(\tilde
{T}_d+\lambda_0 \tilde{F}_d|
\tilde{F}_d\tau_0|\tilde{Z}\bigr)'\proj(
\tilde{U}) \bigl(\tilde {T}_d+\lambda_0 \tilde{F}_d|
\tilde{F}_d\tau_0|\tilde{Z}\bigr)
\nonumber
\\[-8pt]
\\[-8pt]
\nonumber
&=&n\sum_{s}p_s\bigl(T_s+
\lambda_0 F_s|F_s\tau_0|I_p\bigr)'
\tilde {B}\bigl(T_s+\lambda _0 F_s|F_s
\tau_0|I_p\bigr).
\end{eqnarray}
Define $d^*$ to be the design such that
\[
P_{d^*}=\frac{1}{t!}\sum_{\sigma}P_{\sigma d}.
\]
It is easy to show that $d^*$ is a symmetric design and
%
%
\begin{equation}
\label{eqn:11212} C_{d^*,\tau_0,\lambda_0}(\tau,\lambda,\alpha)=\frac{1}{t!}\sum
_{\sigma
}C_{\sigma d,\tau_0,\lambda_0}(\tau,\lambda,\alpha),
\end{equation}
in view of (\ref{eqn:1119}). By Lemma~3.1 of \citet{Kus97} and (\ref
{eqn:11212}), we have
%
%
\begin{equation}
\label{eqn:11213} \frac{1}{t!}\sum_{\sigma}C_{\sigma d,\tau_0,\lambda_0}(
\tau)\leq C_{d^*,\tau_0,\lambda_0}(\tau).
\end{equation}
By (\ref{eqn:11214}) and (\ref{eqn:11213}), we have
\begin{eqnarray*}
\phi_{\delta_{\tau_0},\lambda_0}(d)&=&\frac{1}{t!}\sum_{\sigma
}
\phi _{\delta_{\tau_0},\lambda_0}(\sigma d)
\\
&=&\frac{1}{(t!)^2}\sum_{\sigma}\sum
_{\tilde{\sigma}}\Phi \bigl(C_{\sigma
d,\tilde{\sigma}\tau_0,\lambda_0}(\tau)\bigr)
\\
&=&\frac{1}{(t!)^2}\sum_{\tilde{\sigma}}\sum
_{\sigma}\Phi \bigl(C_{\sigma
d,\tilde{\sigma}\tau_0,\lambda_0}(\tau)\bigr)
\\
&\leq&\frac{1}{t!}\sum_{\tilde{\sigma}}\Phi
\bigl(C_{d^*,\tilde{\sigma
}\tau
_0,\lambda_0}(\tau)\bigr)
\\
&=&\phi_{\delta_{\tau_0},\lambda_0}\bigl(d^*\bigr).
\end{eqnarray*}
\upqed\end{pf}

\begin{remark}
In proving Theorem~\ref{thm:1121} we use the same approach in the proof
of Theorem~3.2 of \citet{Kus97} to derive (\ref{eqn:11213}). However,
the proof of the latter theorem is not rigorous since (3.6) therein
does not hold in general. Actually the gap can be filled by using (\ref
{eqn:11212}) in replacement of (3.6) therein.
\end{remark}

\begin{corol}\label{col:1204}
In approximate design theory, given any value the number $\lambda_0$
and the prior distribution $g$ of $\tau_0$ as long as the latter is
exchangeable, for any design $d$ there exists a symmetric design, say
$d^*$, such that
\[
\phi_{g,\lambda_0}(d)\leq\phi_{g,\lambda_0}\bigl(d^*\bigr).
\]
\end{corol}
\begin{pf}
It is enough to notice that inequality (\ref{eqn:12313}) holds for any
$\tau_0$.
\end{pf}

By Corollary~\ref{col:1204}, there always exists a symmetric design
which is optimal among $\Omega_{p,t,n}$. We define a design $d$ to be
\textit{pseudo symmetric} if all treatments in $d$ are equally replicated
on each period and $C_{dij},1\leq i,j\leq2$ are completely symmetric.
A symmetric design is pseudo symmetric and thus an optimal design in
the subclass of pseudo symmetric designs is automatically optimal among
$\Omega_{p,t,n}$.
%
\begin{proposition}\label{prop:1121}
Regardless the value of $\tau_0$, Fisher's information matrix
$C_{d,\tau_0,\lambda_0}(\tau)$ of a symmetric design $d$ has
eigenvalues of $0,(t-1)^{-1}(c_{d11}-c_{d12}^2/c_{d22})$ and
$(t-1)^{-1}(c_{d11}+2\lambda_0c_{d12}+\lambda_0^2c_{d22})$ with
multiplicities $1,1$ and $t-2$, respectively.
\end{proposition}
\begin{pf}
For a symmetric design $d$, we have $\sum_{s\in{\cal
S}}p_sT_sB_t=0=\break \sum_{s\in{\cal S}}p_sF_sB_t$ and hence
%
%
\begin{equation}
\label{eqn:12142} C_{dij}=\hat{C}_{dij},\qquad 1\leq i,j\leq2,
\end{equation}
in view of (\ref{eqn:1222}). Moreover, these matrices are all
completely symmetric and have row and column sum as zero, which
together with (\ref{eqn:12142}) yields $C_{dij}=c_{dij}B_t/(t-1)$. Due
to $1'\tau_0=0$ and hence $B_t\tau_0=\tau_0$, we have
\begin{eqnarray*}
(t-1)C_{d,\tau_0,\lambda_0}(\tau)=\bigl(c_{d11}+2\lambda
_0c_{d12}+\lambda _0^2c_{d22}
\bigr)B_t-\frac{(c_{d12}+\lambda_0c_{d22})^2}{c_{d22}}\frac
{\tau
_0\tau_0'}{\tau_0'\tau_0}.
\end{eqnarray*}
Let $\{x_1,\ldots,x_{t-2}\}$ be the orthogonal basis which is orthogonal
to both $1$ and $\tau_0$. Then $\{x_1,\ldots,x_{t-2},\tau_0,1\}$ forms the
eigenvectors for the above matrix. Hence, the lemma is concluded.
\end{pf}

\begin{remark}\label{rem:330}
Since $c_{sij}$ is the same for sequences in the same symmetric block
$\langle s\rangle$, we have $c_{dij}=\sum_{s\langle\in\rangle{\cal
S}}p_{\langle s\rangle}c_{sij}$. In view of Corollary~\ref{col:1204}
and Proposition~\ref{prop:1121}, one can derive an optimal design in
two steps. First, we find the optimum value of $p_{\langle s\rangle}$
for all distinct symmetric blocks. Within each symmetric block with
positive $p_{\langle s\rangle}$, we construct a pseudo symmetric
design, and then assemble these designs according to the desired value
of $p_{\langle s\rangle}$. For step one, see equivalence theorems in
Section~\ref{sec:eqthm}. For step two, one can utilize some combinatory
structures such as type I orthogonal arrays, for the latter see Design
6 of \citet{BaiKun06}, for example. For E-criterion, more
general optimal designs could be derived. See Section~\ref{sec:3.3}
for details.
\end{remark}

\begin{remark}
The application of Corollary~\ref{col:1204} and Proposition~\ref
{prop:1121} for A-criterion leads to Proposition~1 of \citet{BaiKun06}.
\end{remark}

\subsection{E-optimality}\label{sec:3.3}
Let ${\cal E}_{g,\lambda_0}(d)$ be the criterion $\phi_{g,\lambda
_0}(d)$ when $\Phi$ therein is evaluated by the second smallest
eigenvalue of the information matrix. We call a design to be ${\cal
E}_{g,\lambda_0}$-optimal if it maximizes ${\cal E}_{g,\lambda_0}(d)$.

\begin{proposition}\label{prop:1204}
In approximate design theory, regardless the value of $\lambda_0$ and
the prior distribution $g$ as long as the latter is exchangeable, a
design $d$ is ${\cal E}_{g,\lambda_0}$-optimal if and only if
${\cal E}_{g,\lambda_0}(d)=ny^*/(t-1)$ with $y^*$ as defined right
before Theorem~\ref{thm:1222}.
\end{proposition}
\begin{pf}
First, it is easy to verify that
\[
c_{d11}-c_{d12}^2/c_{d22}
\leq c_{d11}+2\lambda_0c_{d12}+
\lambda_0^2c_{d22}
\]
for any $\lambda_0$. By Theorem~4.5 of \citet{Kus97}, we have
\begin{eqnarray*}
y^*&=&\min_{-\infty<x<\infty}\sum_{s\in{\cal S}}p_sq_s(x)
\\
&=&n^{-1}\max_d\bigl(c_{d11}-c_{d12}^2/c_{d22}
\bigr).
\end{eqnarray*}
Hence, the proposition is proved in view of Corollary~\ref{col:1204}
and Proposition~\ref{prop:1121}.
\end{pf}

\begin{theorem}\label{thm:1204}
In approximate design theory, regardless of the value of $\lambda_0$ and
the prior distribution $g$ as long as the latter is exchangeable, a
design is ${\cal E}_{g,\lambda_0}$-optimal if there exists a real
number $x$ such that
%
%
\begin{eqnarray}
\sum_{s\in{\cal Q}}p_s[\hat{C}_{s11}+x
\hat{C}_{s12}]&=&\frac
{y^*}{t-1}B_t,\label{eqn:11}
\\
\sum_{s\in{\cal Q}}p_s[\hat{C}_{s21}+x
\hat{C}_{s22}]&=&0,\label {eqn:22}
\\
\sum_{s\in{\cal Q}}p_s\tilde{B}(T_s+xF_s)B_t&=&0,
\label{eqn:33}
\\
\sum_{s\in{\cal Q}}p_s&=&1,
\\
p_s&=&0\qquad \mbox{if }s\notin{\cal Q}.\label{eqn:44}
\end{eqnarray}
\end{theorem}

\begin{pf}
Since $C_{d,\tau_0,\lambda_0}(\tau)$ have column and row sums as zero,
we have
%
%
\begin{equation}
\label{eqn:12044} {\cal E}_{g,\lambda_0}(d)=\E_g \Bigl[\min
_{\ell'1_t=0, \ell'\ell
=1}\ell 'C_{d,\tau_0,\lambda_0}(\tau)\ell \Bigr].
\end{equation}
For a design satisfying (\ref{eqn:11})--(\ref{eqn:44}) we have
%
%
\begin{eqnarray}
C_{d11}+xC_{d12}&=&\frac{ny^*}{t-1}B_t,
\label{eqn:12042}
\\
C_{d21}+xC_{d22}&=&0,\label{eqn:12043}
\end{eqnarray}
in view of (\ref{eqn:1222}). Since $C_{d22}$ is symmetric, (\ref
{eqn:12043}) implies the symmetry of $C_{d21}$ and hence
$C_{12}=C_{21}$. Then by (\ref{eqn:1204}), (\ref{eqn:12042}) and
(\ref
{eqn:12043}), we have
%
%
\begin{equation}
\label{eqn:12045} C_{d,\tau_0,\lambda_0}(\tau)=\frac{ny^*}{t-1}B_t+(
\lambda _0-x)^2C_{d22}-\frac{(\lambda_0-x)^2}{\tau_0'C_{d22}\tau
_0}C_{d22}
\tau _0\tau_0'C_{d22}.
\end{equation}
Let $\{0,a_1,\ldots,a_{t-1}\}$ be the eigenvalues of $C_{d22}$ with
corresponding normalized eigenvectors $\{1_t,\ell_1,\ldots,\ell_{t-1}\}$,
then we have $C_{d22}=\sum^{t-1}_{i=1}a_i\ell_i\ell_i'$. Since $\tau
_0'1_t=0$, we have the representation $\tau_0=\sum^{t-1}_{i=1}c_i\ell
_i$. For any vector $\ell$ with $\ell'1_t=0$ and $\ell'\ell=1$, we have
the expression of $\ell=\sum^{t-1}_{i=1}b_i\ell_i$ with the restriction
$\sum^{t-1}_{i=1}b_i^2=1$, the equation $\ell'B_t\ell=1$, and hence by
(\ref{eqn:12045})
%
%
\begin{eqnarray}
\label{eqn:12046} \ell'C_{d,\tau_0,\lambda_0}(\tau)\ell&=&
\frac{ny^*}{t-1}+(\lambda _0-x)^2\sum
^{t-1}_{i=1}a_ib_i^2-
\frac{(\lambda_0-x)^2}{\sum^{t-1}_{i=1}a_ic_i^2} \Biggl(\sum^{t-1}_{i=1}a_ib_ic_i
\Biggr)^2
\nonumber
\\
&=&\frac{ny^*}{t-1}+\frac{(\lambda_0-x)^2}{\sum^{t-1}_{i=1}a_ic_i^2} \Biggl[ \Biggl(\sum
^{t-1}_{i=1}a_ib_i^2
\Biggr) \Biggl(\sum^{t-1}_{i=1}a_ic_i^2
\Biggr)- \Biggl(\sum^{t-1}_{i=1}a_ib_ic_i
\Biggr)^2 \Biggr]
\\
&\geq&\frac{ny^*}{t-1},\nonumber
\end{eqnarray}
the equality holds if and only if $\ell=\tau_0/\Vert \tau_0\Vert $. The
theorem is concluded in view of Proposition~\ref{prop:1204}, (\ref
{eqn:12044}) and (\ref{eqn:12046}).
\end{pf}

\begin{remark}
The advantage of Theorem~\ref{thm:1204} is that the design therein is
optimal for any $\lambda_0$ while the A-optimality of totally balanced
design [\citet{BaiKun06}] requires the condition of $-1\leq
\lambda_0\leq\lambda^*$.
\end{remark}
As a direct result of Theorems \ref{thm:1222} and \ref{thm:1204}, we have
the following corollary.
%
\begin{corol}\label{cor:12182}
In approximate design theory, regardless the value of $\lambda_0$ and
the prior distribution $g$ as long as the latter is exchangeable, a
universally optimal design for model (\ref{model3}) is also ${\cal
E}_{g,\lambda_0}$-optimal for model (\ref{model4}).
\end{corol}

\begin{theorem}
The variable $x$ in (\ref{eqn:11})--(\ref{eqn:44}) takes the unique
value of $x^*$, which is defined right above Theorem~\ref{thm:1222}.
\end{theorem}

\begin{pf}
Given (\ref{eqn:11})--(\ref{eqn:44}), we have by (\ref{eqn:12042}) and
(\ref{eqn:12043}) that
\begin{eqnarray*}
C_{d11}-C_{d12}(C_{d22})^-C_{d21}&=&C_{d11}+xC_{d12}(C_{d22})^-C_{d22}
\\
&=&C_{d11}+xC_{d12}
\\
&=&\frac{ny^*}{t-1}B_t,
\end{eqnarray*}
which indicates that $d$ is universally optimal for model (\ref
{model3}) in view of Theorem~\ref{thm:1222}. Hence, we have $x=x^*$ by
Theorem~\ref{thm:1222}.
\end{pf}

As a direct result of Theorem~\ref{thm:12222}, Corollary~\ref
{cor:12182} and Remark~\ref{rem:330}, we have
Corollary~\ref{co3}.
%
\begin{corol}\label{co3}
In approximate design theory, regardless of the value of $\lambda_0$ and
the prior distribution $g$ as long as the latter is exchangeable, a
pseudo symmetric design is ${\cal E}_{g,\lambda_0}$-optimal if it
satisfies (\ref{eqn:12222})--(\ref{eqn:12223}).
\end{corol}

\subsection{Equivalence theorems}\label{sec:eqthm}

In order to introduce the following results, we define
$x_d=-c_{d12}/c_{d22}$ and $q_d(x)=\sum_{s\in{\cal S}}p_sq_s(x)$. Then
we have
\begin{eqnarray*}
nq_d(x_d)&=&c_{d11}-c_{d12}^2/c_{d22},
\\
nq_d(\lambda_0)&=&c_{d11}+2
\lambda_0 c_{d12}+\lambda_0^2c_{d22}.
\end{eqnarray*}
For a $t\times t$ matrix $C$ with eigenvalues $0=a_0\leq a_1\leq
a_2\leq\cdots\leq a_{t-1}$, define the criterion functions
\begin{eqnarray*}
\Phi_A(C)&=&(t-1) \Biggl(\sum^{t-1}_{i=1}a_i^{-1}
\Biggr)^{-1},
\\
\Phi_D(C)&=& \Biggl(\prod_{i=1}^{t-1}a_i
\Biggr)^{1/(t-1)},
\\
\Phi_T(C)&=&(t-1)^{-1}\sum^{t-1}_{i=1}a_i.
\end{eqnarray*}
Let ${\cal A}_{g,\lambda_0}(d)$, ${\cal D}_{g,\lambda_0}(d)$ and
${\cal
T}_{g,\lambda_0}(d)$ be the criterion $\phi_{g,\lambda_0}(d)$ when
$\Phi
$ therein is evaluated by $\Phi_A$, $\Phi_D$ and $\Phi_T$,
respectively. We call a design to be ${\cal A}_{g,\lambda_0}$-optimal
if it maximizes ${\cal A}_{g,\lambda_0}(d)$. Definitions for optimality
of ${\cal D}_{g,\lambda_0}$ and ${\cal T}_{g,\lambda_0}$ are similar.

\begin{theorem}\label{thm:D}
In approximate design theory, regardless of the value of $\lambda_0$ and
the prior distribution $g$ as long as the latter is exchangeable, a
pseudo symmetric design $d$ is ${\cal D}_{g,\lambda_0}$-optimal if and
only if
%
%
\begin{eqnarray}
\label{eqn:330} \max_{s\in{\cal S}} \biggl(\frac{1}{t-1}
\frac
{q_s(x_d)}{q_d(x_d)}+\frac
{t-2}{t-1}\frac{q_s(\lambda_0)}{q_d(\lambda_0)} \biggr)&=&1.
\end{eqnarray}
Moreover, the sequences in design $d$ attain the maximum in (\ref{eqn:330}).
\end{theorem}
\begin{pf}
For a real number $x$, let $\eta(\xi_1,\xi_2,x)=q_{\xi_2}(x)/q_{\xi
_1}(x)$ with $\xi_1$ and $\xi_2$ being either a design or a sequence.
Also define $\psi
(P_d)=\operatorname{log}((c_{d11}-c_{d12}^2/c_{d22}))+(t-2)\operatorname{log}(c_{d11}+2\lambda
_0c_{d12}+\lambda_0^2c_{d22})$. By Theorem~\ref{thm:12222}, Corollary~\ref{cor:12182}, Remark~\ref{rem:330} and the concavity of D-criterion,
a pseudo symmetric design $d^*$ is ${\cal D}_{g,\lambda_0}$-optimal if
and only if for any other design $d$ we have
%
%
\begin{eqnarray}
0&\geq&\lim_{\delta\rightarrow0}\frac{\psi((1-\delta
)P_{d^*}+\delta
P_{d})-\psi(P_{d^*})}{\delta}
\nonumber
\\[-8pt]
\\[-8pt]
\nonumber
&=&\eta \bigl(d^*,d,x_{d^*} \bigr)+(t-2)\eta\bigl(d^*,d,\lambda
_0\bigr)-(t-1).\label{eqn:3301}
\end{eqnarray}
Take $d$ in (\ref{eqn:3301}) to be a design consist of a single
sequence $s$, we have
%
%
\begin{equation}
\label{eqn:3304} \max_{s\in{\cal S}} \biggl(\frac{1}{t-1}\eta
\bigl(d^*,s,x_{d^*} \bigr)+\frac{t-2}{t-1}\eta\bigl(d^*,s,
\lambda_0\bigr) \biggr)\leq 1.
\end{equation}
Observe that
%
%
\begin{equation}
\label{eqn:3302} \eta (d,d,x_{d} )=1=\eta(d,d,\lambda_0),
\end{equation}
then we have
%
%
\begin{equation}
\label{eqn:3303} \max_{s\in{\cal S}} \biggl(\frac{1}{t-1}\eta
(d,s,x_{d} )+\frac
{t-2}{t-1}\eta(d,s,\lambda_0)
\biggr)\geq 1.
\end{equation}
The theorem is completed in view of (\ref{eqn:3304}), (\ref{eqn:3302})
and (\ref{eqn:3303}).
\end{pf}

\begin{theorem}\label{thm:401}
In approximate design theory, regardless of the value of $\lambda_0$ and
the prior distribution $g$ as long as the latter is exchangeable, a
pseudo symmetric design $d$ is ${\cal A}_{g,\lambda_0}$-optimal if and
only if
%
%
\begin{equation}
\label{eqn:3306} \max_{s\in{\cal S}} \biggl(\pi_d
\frac{q_s(x_d)}{q_d(x_d)}+(1-\pi _d)\frac
{q_s(\lambda_0)}{q_d(\lambda_0)} \biggr)=1,
\end{equation}
where $\pi_d=q_d(\lambda_0)/(q_d(\lambda_0)+(t-2)q_d(x_d))$. Moreover,
the sequences in design $d$ attain the maximum in (\ref{eqn:3306}).
\end{theorem}
%
\begin{remark}
Theorem~\ref{thm:401} is essentially a generalization of the result of
\citet{BaiKun06}.
\end{remark}

\begin{theorem}\label{thm:E}
In approximate design theory, regardless of the value of $\lambda_0$ and
the prior distribution $g$ as long as the latter is exchangeable, a
pseudo symmetric design $d$ is ${\cal E}_{g,\lambda_0}$-optimal if and
only if
%
%
\begin{equation}
\label{eqn:401} \max_{s\in{\cal S}}\frac{q_s(x_d)}{q_d(x_d)}=1.
\end{equation}
Moreover, the sequences in design $d$ attain the maximum in (\ref{eqn:401}).
\end{theorem}
%
\begin{remark}
In fact, (\ref{eqn:401}) is equivalent to (\ref{eqn:12222})--(\ref
{eqn:12223}).
\end{remark}

\begin{theorem}\label{thm:T}
In approximate design theory, regardless of the value of $\lambda_0$ and
the prior distribution $g$ as long as the latter is exchangeable, a
pseudo symmetric design $d$ is ${\cal T}_{g,\lambda_0}$-optimal if and
only if
%
%
\begin{equation}
\label{eqn:3305} \max_{s\in{\cal S}}\frac{q_s(x_d)+(t-2)q_s(\lambda
_0)}{q_d(x_d)+(t-2)q_d(\lambda_0)}=1.
\end{equation}
Moreover, the sequences in design $d$ attain the maximum in (\ref{eqn:3305}).
\end{theorem}

\subsection{\texorpdfstring{Estimation of $\lambda_0$}{Estimation of lambda 0}}\label{sec:3.5}

By the Cram\'er--Rao inequality, the variance of an unbiased estimator of
$\lambda_0$ is bounded by the reciprocal of
%
%
\begin{equation}
\label{eqn:0102} C_{d,\tau_0,\lambda_0}(\lambda)=\tau_0'
\tilde{F}_d'\proj\bigl(\tilde {T}_d+
\lambda_0 \tilde{F}_d|\tilde{Z}|\tilde{U}\bigr)
\tilde{F}_d\tau_0,
\end{equation}
achievable by MLE asymptotically. Define $A_{d11}=F_d'(B_n\otimes
\tilde
{B})F_d$, $A_{d11}=F_d'(B_n\otimes\tilde{B})(T_d+\lambda_0 F_d)$,
$A_{d21}=A_{d12}'$ and $A_{d22}=(T_d+\lambda_0 F_d)'(B_n\otimes\tilde
{B})(T_d+\lambda_0 F_d)$. Straightforward calculations show that
$C_{d,\tau_0,\lambda_0}(\lambda)=\tau_0'A_d\tau_0$ where
$A_d=A_{d11}-A_{d12}(A_{d22})^-A_{21}$.

As in (\ref{eqn:0103}) we define $\varphi_{g,\lambda_0}(d)=\E_g\tau
_0'A_d\tau_0$, where the expectation is taken with respect to the prior
distribution measure $g$. Then we have $\varphi_{\delta_{\tau
_0},\lambda
_0}(d)=\tau_0'\bar{A}_d\tau_0$ where $\bar{A}_d=\frac{1}{t!}\sum_{\sigma
}S_{\sigma}'A_dS_{\sigma}$.

For each sequence $s$, define $\hat{A}_{sij}=B_tG_i'\tilde{B}G_jB_t,
1\leq i,j\leq2$ with $G_1=F_s$ and $G_2=T_s+\lambda_0F_s$. By direct
calculations, we have
%
%
\begin{equation}
A_{dij}=\hat{A}_{dij}-nG_i'
\tilde{B}G_j,\qquad 1\leq i,j\leq2,\label{eqn:1222}
\end{equation}
where $\hat{A}_{dij}=n\sum_{s\in{\cal S}}p_s\hat{A}_{sij},1\leq
i,j\leq2$, $G_1=\sum_{s\in{\cal S}}p_sF_sB_t$\vspace*{2pt} and $G_2=\break \sum_{s\in
{\cal S}}p_s(T_s+\lambda_0F_s)B_t$.

Further define $h_{sij}=\operatorname{tr}(\hat{A}_{sij})$, $h_{dij}=\operatorname{tr}(\hat
{A}_{dij})=n\sum_{s\in{\cal S}}p_sh_{sij}$, the quadratic function
$r_s(x)=h_{s11}+2h_{s12}x+h_{s22}x^2$, $r(x)=\max_{s}r_s(x)$,
$y_0=\break \min_{-\infty<x<\infty}r(x)$, $x_0$ to be the unique solution of $r(x)=y_0$
and ${\cal R}=\{s\in{\cal S}|r_s(x_0)=y_0\}$. Now we have the
following theorem.

\begin{theorem}\label{thm:12223}
Given any $-\infty<\lambda_0<\infty$, a design maximizes $\varphi
_{\delta_{\tau_0},\lambda_0}(d)$ for any $\tau_0$ if
%
%
\begin{eqnarray}
\sum_{s\in{\cal R}}p_s[\hat{A}_{s11}+x_0
\hat{A}_{s12}]&=&\frac
{y_0}{t-1}B_t,\label{eqn:1}
\\
\sum_{s\in{\cal R}}p_s[\hat{A}_{s21}+x_0
\hat{A}_{s22}]&=&0,\label {eqn:2}
\\
\sum_{s\in{\cal R}}p_s\tilde{B}\bigl[
\hat{F}_s+x_0(\hat{T}_s+
\lambda_0 \hat {F}_s)\bigr]&=&0,\label{eqn:3}
\\
\sum_{s\in{\cal R}}p_s&=&1,
\\
p_s&=&0\qquad\mbox{if }s\notin{\cal R}.\label{eqn:4}
\end{eqnarray}
\end{theorem}

\begin{pf}
Note that for any design $d$, there exists a symmetric design $d^*$
with $\bar{A}_d\leq A_{d^*}$ by the same argument as for (\ref
{eqn:11213}). For a symmetric design, we have
$A_d=(t-1)^{-1}(h_{d11}-h_{d12}^2/h_{d22})B_t$. By similar arguments as
in proof of Theorem~4.4 of \citet{Kus97}, we have $\max_d(h_{d11}-h_{d12}^2/h_{d22})=ny_0$.
By direct calculations, we know
that (\ref{eqn:1})--(\ref{eqn:4}) implies $A_d=ny_0B_t/(t-1)$ and hence
the theorem is proved.
\end{pf}

\begin{corol}
Given any value of the real number $\lambda_0$, a design maximizes
$\varphi_{g,\lambda_0}(d)$ for any exchangeable prior distribution $g$
if it satisfies \mbox{(\ref{eqn:1})--(\ref{eqn:4})}.
\end{corol}
\begin{pf}
The necessity is immediate. For sufficiency, it is enough to note that
the design satisfying (\ref{eqn:1})--(\ref{eqn:4}) does not depend on
$\tau_0$.
\end{pf}

\begin{corol}
Given any value of the real number $\lambda_0$, a design maximizes
$\varphi_{g,\lambda_0}(d)$ for any exchangeable prior distribution $g$
if it is a symmetric design with
%
%
\begin{eqnarray}
\sum_{s\langle\in\rangle{\cal R}}p_{\langle s\rangle
}r_s'(x_0)&=&0,
\label{eqn:1224}
\\
\sum_{s\langle\in\rangle{\cal R}}p_{\langle s\rangle}&=&1,
\\
p_s&=&0\qquad\mbox{if }s\notin{\cal R},\label{eqn:12242}
\end{eqnarray}
where $r_s'(x)$ is the derivative of $r_s(x)$ with respect to $x$.
\end{corol}

\begin{pf}
It is enough to show that (\ref{eqn:1224})--(\ref{eqn:12242}) implies
(\ref{eqn:1})--(\ref{eqn:4}). The proof of the latter is analogous to
that of Theorem~\ref{thm:12222}.
\end{pf}

A general necessary and sufficient optimality condition is given by
the following.
%
\begin{theorem}
Given any value of the real number $\lambda_0$, a design maximizes
$\varphi_{g,\lambda_0}(d)$ for any exchangeable prior distribution $g$
if and only if
$\operatorname{tr}(A_d)=ny_0$.
\end{theorem}
\begin{pf}
Note that $\bar{A}_d=\operatorname{tr}(A_d)B_t/(t-1)$, hence we have $\varphi
_{\delta
_{\tau_0},\lambda_0}(d)=\tau_0'\tau_0\operatorname{tr}(A_d)/(t-1)$ and hence
%
%
\begin{equation}
\label{eqn:12314} \varphi_{g,\lambda_0}(d)=\frac{\E_g(\tau_0'\tau_0)}{t-1}\operatorname{tr}(A_d).
\end{equation}
Through the proof of Theorem~\ref{thm:12223}, we know that $\max_d
\operatorname{tr}(A_d)=ny_0$, which together with (\ref{eqn:12314}) proves the theorem.
\end{pf}

\section{Examples}\label{sec:exm}
In the spirit of Theorems \ref{thm:D}--\ref{thm:T} and Remark~\ref
{rem:330}, we consider examples of optimal designs in the format of
pseudo symmetric designs, even though a more general format could be
proposed for E-optimality due to Theorem~\ref{thm:1204}. Let $m$ to be
the total number of distinct symmetric blocks and suppose $s_1,s_2,\ldots,s_m$
are the representative sequences for each of the symmetric blocks. For
a design $d$, define the vector $P_{\langle d\rangle}=(p_{\langle
s_1\rangle},p_{\langle s_2\rangle},\ldots,p_{\langle s_m\rangle})$. Then
two pseudo symmetric designs with the same $P_{\langle d\rangle}$ will
have the same $\phi_{g,\lambda_0}(d)$ for any $\phi,g,\lambda_0$ as
long as $g$ is exchangeable. In particular, they are equivalent in
terms of ${\cal A}_{g,\lambda_0}$-, ${\cal D}_{g,\lambda_0}$-, ${\cal
E}_{g,\lambda_0}$- and ${\cal T}_{g,\lambda_0}$-optimality. In the
sequel, we will mainly focus on the determination of $P_{\langle
d\rangle}$ based on the equivalence theorems \ref{thm:D}--\ref{thm:T}.
A general algorithm could be found in the supplemental article [Zheng (\citeyear{supp})].

For the following examples, we consider first order autocorrelation for
within subject covariance matrix, namely $\Sigma=(\rho\I
_{|i-j|=1}+\I
_{i=j})_{1\leq,i,j\leq p}$, where $\I$ is the indicator function.
Hence, $\rho=0$ implies $\Sigma=I_p$. Following \citet{Kus98}, we
define two special symmetric blocks. The symmetry block $\langle
\mathrm{di}\rangle$ consists of all sequences having distinct treatments in the
$p$ periods. The symmetry block $\langle \mathrm{re}\rangle$ consists of all
sequences having distinct treatments in the first $p-1$ periods, with
the treatment in period $p-1$ repeating in period $p$. All examples
given below are pseudo symmetric designs except otherwise specified.
For ease of illustration by examples, we only consider $(\rho,\lambda
_0)\in\{-1/2,0,1/2\}\times[-1,1]$, even though other values of $(\rho
,\lambda_0)$ does not cause extra difficulty. Throughout this section,
$g$ is exchangeable unless otherwise specified.

\textit{Case of $(p,t)=(3,3)$}:
Let $d_1$ be a design with $p_{\langle \mathrm{re}\rangle}=1/6$ and $p_{\langle
\mathrm{di}\rangle}=5/6$. See, for instance, Example~1 of \citet{Kus98} with
$n=36$ subjects. Define $d_2$ to be a design with $p_{\langle \mathrm{di}\rangle
}=1$, which requires $n$ to be a multiple of $6$ as an exact design.
When $\rho=0$, Theorem~\ref{thm:E} shows the ${\cal E}_{g,\lambda
_0}$-optimality of $d_1$ for any $\lambda_0$ and \citet{BaiKun06} shows the ${\cal A}_{g,\lambda_0}$-optimality of $d_2$ when
$-1\leq\lambda_0\leq\lambda^*=0.34375$. In fact, one can verify by
Theorems \ref{thm:401}, \ref{thm:D} and \ref{thm:T} that $d_2$ is even
${\cal A}_{g,\lambda_0}$-, ${\cal D}_{g,\lambda_0}$- and ${\cal
T}_{g,\lambda_0}$-optimal when $-1\leq\lambda_0\leq0.394$. At
$\lambda
_0=0.5$, $d_1$ is optimal under all four criteria. When we tune $\rho$
to be $1/2$, $d_2$ is optimal under all four criteria for $-0.75\leq
\lambda_0\leq1$. When we tune $\rho$ to be $-1/2$, $d_2$ is still
${\cal A}_{g,\lambda_0}$-, ${\cal D}_{g,\lambda_0}$- and ${\cal
T}_{g,\lambda_0}$-optimal for small and negative values of $\lambda_0$,
while the design with $p_{\langle \mathrm{re}\rangle}=2/9$ and $p_{\langle
\mathrm{di}\rangle}=7/9$ is ${\cal E}_{g,\lambda_0}$-optimal. For moderate
positive $\lambda_0$, designs for four criteria are all different, but
they all consists of small portion of ${\langle \mathrm{re}\rangle}$ and large
portion of ${\langle \mathrm{di}\rangle}$. All these designs are highly
efficient for all criteria; see Table~\ref{table33}, for example.

\begin{table}
\tablewidth=250pt
\caption{Efficiency of $d_1$ and $d_2$ under ${\cal A}_{g,\lambda_0}$,
${\cal D}_{g,\lambda_0}$, ${\cal E}_{g,\lambda_0}$, and ${\cal
T}_{g,\lambda_0}$-criteria for the case of $(p,t)=(3,3)$ when $\rho
=\lambda_0=0$}\label{table33}
\begin{tabular*}{250pt}{@{\extracolsep{\fill}}lcccc@{}}
\hline
\textbf{Design} &\textbf{A}&\textbf{D}&\textbf{E}&\textbf{T}\\
\hline
$d_1$& 0.9782& 0.9752& 1\phantom{0000.} &0.9722\\
$d_2$& 1\phantom{0000.}& 1\phantom{0000.} &0.9931& 1\phantom{0000.}\\
\hline
\end{tabular*}
\end{table}

Without surprise, $\phi_{g,\lambda_0}$-optimal design for exchangeable
$g$ is not necessarily optimal when $g$ is not exchangeable. We
consider the prior distribution of $g=g_1$ which puts all its mass on
the single point $\tau_0=(0,1,-1)'$. When $n=36$, derive $d_{1'}$ from
$d_1$ by replacing one sequence of $123$ therein by $323$, it turns out
that $d_{1'}$ is $1.66\%$ more ${\cal E}_{g_1,0}$-efficient when
$\lambda_0=\rho=0$. However, in practice, one does not have accurate
information of $\tau_0$. Exchangeable prior distribution of $\tau_0$
actually accounts for the case when nothing is known about $\tau_0$. A
further justification is that symmetry is usually a nice feature. If we
search among pseudo symmetric designs, the designs as proposed in this
paper would be optimal under the corresponding criterion for any
arbitrary $g$, which is not necessarily exchangeable. To see this, note
that $C_{dij}$'s are all completely symmetric for a pseudo symmetric
design. Hence, it is easily seen, by examining the proof of Proposition~\ref{prop:1121},
that the value of $\phi_{g,\lambda_0}(d)$ is
independent of the distribution $g$ for both $d_1$ and $d_2$ regardless
the value of $\lambda_0$ as well as the criterion function $\Phi$.

\textit{Case of $(p,t)=(3,4)$}:
Let $d_3$ be a design with $p_{\langle \mathrm{re}\rangle}=1/8$ and $p_{\langle
\mathrm{di}\rangle}=7/8$. Define $d_4$ to be a design with $p_{\langle
\mathrm{di}\rangle
}=1$, which requires $n$ to be a multiple of $12$ as an exact design.
When $\rho=0$, Theorem~\ref{thm:E} shows the ${\cal E}_{g,\lambda
_0}$-optimality of $d_3$ for any $\lambda_0$ and \citet{BaiKun06} shows the ${\cal A}_{g,\lambda_0}$-optimality of $d_4$ when
$-1\leq\lambda_0\leq\lambda^*=0.4455$. In fact, one can verify by
Theorems \ref{thm:D}, \ref{thm:401} and \ref{thm:T} that $d_4$ is even
${\cal A}_{g,\lambda_0}$-, ${\cal D}_{g,\lambda_0}$- and ${\cal
T}_{g,\lambda_0}$-optimal when $-1\leq\lambda_0\leq0.463$. Similarly
at $\lambda_0=0.5$, $d_3$ is optimal under all four criteria. When we
tune $\rho$ to be $1/2$, $d_4$ is optimal under all four criteria for
$-0.35\leq\lambda_0\leq1$. It is still ${\cal A}_{g,\lambda_0}$-,
${\cal D}_{g,\lambda_0}$- and ${\cal E}_{g,\lambda_0}$-optimal and
highly ${\cal T}_{g,\lambda_0}$-efficient for $-1\leq\lambda_0<-0.35$.
When $\rho=-0.5$, similar phonomania as for the case of $(p,t)=(3,4)$
is observed.

\textit{Case of $(p,t)=(3,5)$}:
We have similar observations as for case of $(p,t)=(3,4)$, except that
the portion of ${\langle \mathrm{re}\rangle}$ becomes further smaller. This
trend projects to larger values of $t$.

\textit{Case of $(p,t)=(4,3)$}:
When $\rho=0$, the design with $p_{\langle \mathrm{re}\rangle}=1$ is optimal
under all four criteria for $0\leq\lambda_0\leq1$. For negative
$\lambda_0$, designs are different for different criteria. However,
they typically consist of symmetric blocks of $\langle1232\rangle$
and $\langle \mathrm{re}\rangle$. Table~\ref{table43} shows the performance of
these designs for $\lambda_0=-0.5$. Designs therein are identified by
$p_{\langle \mathrm{re}\rangle}=1-p_{\langle1232\rangle}$. When $\rho$ is
nonzero, symmetric blocks of $\langle1123\rangle$, $\langle
1231\rangle
$,$\langle1232\rangle$ and $\langle \mathrm{re}\rangle$ will appear in
different optimal designs. Note that \citeauthor{BaiKun06}'s (\citeyear{BaiKun06}) result
does not apply to this case since they deal with $3\leq p\leq t$.

\begin{table}
\tablewidth=250pt
\caption{Efficiency of designs under ${\cal A}_{g,\lambda_0}$, ${\cal
D}_{g,\lambda_0}$, ${\cal E}_{g,\lambda_0}$, and ${\cal T}_{g,\lambda
_0}$-criteria for the case of $(p,t)=(4,3)$ when $\rho=0$ and $\lambda
_0=-0.5$}\label{table43}
\begin{tabular*}{250pt}{@{\extracolsep{\fill}}lcccc@{}}
\hline
\multicolumn{1}{@{}l}{$\bolds{p_{\langle1232\rangle}}$}&\textbf{A}&\textbf{D}&\textbf{E}&\textbf{T}\\
\hline
0.4729& 1\phantom{0000.} &0.9964& 0.9553& 0.9768\\
0.6330& 0.9952& 1\phantom{0000.}& 0.9199& 0.9888\\
0\phantom{0000.}& 0.9636& 0.9475& 1\phantom{0000.}& 0.9167\\
1\phantom{0000.}& 0.9422& 0.9785& 0.8000& 1\phantom{0000.}\\
\hline
\end{tabular*}
\end{table}

\textit{Case of $(p,t)=(4,4)$}: When $\rho=0$, the design with
$p_{\langle
\mathrm{re}\rangle}=1/12$ and $p_{\langle \mathrm{di}\rangle}=11/12$ is ${\cal
E}_{g,\lambda_0}$-optimal for all $\lambda_0$, while the design with
$p_{\langle \mathrm{di}\rangle}=1$ is ${\cal A}_{g,\lambda_0}$-, ${\cal
D}_{g,\lambda_0}$- and ${\cal T}_{g,\lambda_0}$-optimal for any
$\lambda
_0$ between $-1$ and $0.318\ (>\lambda^*)$. Interestingly, the design
with $p_{\langle \mathrm{re}\rangle}=1$ is ${\cal A}_{g,\lambda_0}$-, ${\cal
D}_{g,\lambda_0}$- and ${\cal T}_{g,\lambda_0}$-optimal for
$0.625\leq
\lambda_0\leq1$. For $0.318<\lambda<0.625$, the optimal designs
consist of $\langle \mathrm{re}\rangle$ and $\langle \mathrm{di}\rangle$ with the
proportion depending on different criteria. When $\rho=0.5$, the design
with $p_{\langle \mathrm{di}\rangle}=1$ is optimal under all the four criteria
for $0.368\leq\lambda_0\leq1$ and ${\cal E}_{g,\lambda_0}$-optimal
for all $\lambda_0$. When $\rho=-0.5$, the design with $p_{\langle
1123\rangle}=1$ (resp., $p_{\langle \mathrm{re}\rangle}=1$) is ${\cal
E}_{g,\lambda_0}$-optimal for all $\lambda_0$ and also optimal under
the other three criteria for $\lambda_0$ close to zero (reps. $0.3$).
For moderate negative value of $\lambda_0$, the design with
$p_{\langle
\mathrm{di}\rangle}=1$ is optimal under these three
criteria.\looseness=1

\textit{Case of $(p,t)=(4,5)$}: Similar observation as the case of
$(p,t)=(4,4)$ except that the symmetric metric $\langle1122 \rangle$
appears as small proportion in optimal designs when $\rho=-0.5$ and
$\lambda_0$ takes a positive moderate value.

\textit{Case of $(p,t)=(5,3)$}: When $\rho=0$, the design with
$p_{\langle
12233\rangle}=2/5$ and $p_{\langle12332\rangle}=3/5$ is ${\cal
E}_{g,\lambda_0}$-optimal for all $\lambda_0$ and also optimal under
the other three criteria when $\lambda_0$ is in a neighborhood of zero.
For other values of $\rho$ and $\lambda_0$, there is no specific
symmetric block which will dominate, but we observe that all sequences
in the optimal designs contain all three treatments.\looseness=1

\textit{Case of $(p,t)=(6,2)$}: It is well known that $t=2$ indicates the
equivalence of all optimality criteria for the classical model. For
proportional model, this is also true. To see this, Proposition~\ref
{prop:1121} shows that the information matrix $C_{d,\tau_0,\lambda
_0}(\tau)$ only has one positive eigenvalue
$(t-1)^{-1}(c_{d11}-c_{d12}^2/c_{d22})$ with multiplicity 1. Hence, the
optimal design will be irrelevant of optimality criteria as well as the
value of $\lambda_0$. When $\rho=0$, the design with $p_{\langle
111222\rangle}=5/8$ and $p_{\langle121212\rangle}=3/8$ is optimal
under all the four criteria for any $\lambda_0$. An exact design with
$16$ runs is given as\vspace*{6pt}
\[
\begin{array} {cccccccccccccccc} 1&1&1&1&1&2&2&2&2&2&1&1&1&2&2&2
\\
1&1&1&1&1&2&2&2&2&2&2&2&2&1&1&1
\\
1&1&1&1&1&2&2&2&2&2&1&1&1&2&2&2
\\
2&2&2&2&2&1&1&1&1&1&2&2&2&1&1&1
\\
2&2&2&2&2&1&1&1&1&1&1&1&1&2&2&2
\\
2&2&2&2&2&1&1&1&1&1&2&2&2&1&1&1 \end{array} \vspace*{6pt}
\]
When $\rho=0.5$, the design with $p_{\langle122121\rangle}=1$ is
optimal. When $\rho=-0.5$, the design with $p_{\langle111222\rangle
}=2/11$ and $p_{\langle122211\rangle}=9/11$ is optimal. Based on
Corollary~\ref{col:510}, these designs are also universally optimal for
the classical model.

\begin{corol}\label{col:510}
For any $\Sigma$, a pseudo symmetric design which is ${\cal
E}_{g,\lambda_0}$-optimal for the proportional model is also
universally optimal for the classical model.
\end{corol}

\section*{Acknowledgments}
We are grateful to the Associate Editor and two referees
for their constructive comments on earlier versions of this manuscript.

\begin{supplement}[id=suppA]
\stitle{Appendix for
optimal crossover designs for the
proportional model}
\slink[doi]{10.1214/13-AOS1148SUPP} 
\sdatatype{.pdf}
\sfilename{aos1148\_supp.pdf}
\sdescription{This document is to
provide a general algorithm to derive optimal $P_{\langle d \rangle}$ for arbitrary values of
$\lambda_0$ and $\Sigma$ based on the equivalence theorems.}
\end{supplement}


\printaddresses

\end{document}